\documentclass[a4paper,10pt]{article}
\usepackage{amssymb}
\usepackage[T2A]{fontenc}
\usepackage[cp1251]{inputenc}
\usepackage{amsmath}

\usepackage[dvips]{graphicx}
\usepackage[english,russian]{babel}

\newtheorem{theorem}{Теорема}[section]

\newtheorem{definitionhead}[theorem]{Определение}

\begin{document}
\author{А.~Я.~Белов}
\title{Научные конференции школьников}

\maketitle

{
\flushright{Посвящается памяти Р.~В.~Плыкина, организатора и вдохновителя Всероссийской конференции учащихся ``Юность, наука, культура''}

}

%\bigskip

%Важность организации учебно-исследовательской работы учащихся трудно переоценить. Ее значение не исчерпывается задачами чисто професионального
%образования. Самостоятельно полученный результат дает молодому человеку дополнительную внутреннюю опору. Особенностью отечественного образования является работа с ранним возрастом.

\section{Введение }

%\large
Во внешкольной работе по математике распространено два подхода. Большая часть
внеклассной работы школьников прямо или косвенно связана с олимпиадами или иными соревнованиями по решению четко поставленных задач в жестких временных рамках.  Даже если руководитель кружка заявляет, что не занимается подготовкой к олимпиадам, материалы кружка в сильной степени основаны на олимпиадных задачах. Олимпиадный подход имеет свои достоинства и недостатки, ему посвящена статья \cite{dver}.

Второй подход представляет собой проектную деятельность и тесно связан с научными конференциями учащихся. Проектный подход, в отличие от задачного подхода, менее специфичен для математики и выработанные умения легче переносятся и на другие науки. Научные конференции школьников зачастую проводятся сразу по нескольким предметам, математика при этом является
одним из направлений. Они являются альтернативной формой деятельности, позволяющей увидеть и развить в школьнике полезные для научной работы качества, которые олимпиады не раскрывают. Эти два подхода взаимно дополнительны. Надо стремиться к тому, чтобы преподаватели, исповедующие тот или иной подход, не конкурировали между собой, а сотрудничали, понимали друг друга. В выигрыше окажется ученик.

Автор в течение более чем 10 лет был руководителем математической секции конференции ``Поиск'', проходившей в ДНТТМ (Доме научно-технического творчества молодежи), филиале Московского Городского Дворца детского (юношеского) творчества, а также членом экспертного совета Всероссийской конференции учащихся ``Юность, наука, культура'', проходившей по инициативе Р.~В.~Плыкина. Кроме того, автор свыше 20 лет являлся руководителем проектов на Летней конференции турнира городов. Данная статья -- попытка осмыслить, систематизировать и передать опыт работы, который ему не безразличен. Этот опыт относится к конференциям по математике, но многие аспекты являются общими и важными для научных конференций школьников и по другим предметам.

\medskip
{\bf Структура статьи.} После истории вопроса мы обсуждаем академические ценности, для достижения которых особенно полезна проектная деятельность. Затем обсуждаются критические замечания, методика проведения, в особенности стратегические вещи  -- принципы организации жюри и оценки докладов.
В заключение говорится о  перспективах.

\medskip

%\subsection
{\bf История вопроса.} В конце 20-х --- начале 30-х годов
 А.~Н.~Колмогорову (в соответствии со стилем школы Н.~Н.~Лузина) удалось провести реформу профессионального математического образования.
\footnote{Автору представляется, что эта реформа удалась за счет наличия достаточного числа относительно приличных  школ.
  Было бы весьма интересно провести сравнительный анализ знаниий, умений и навыков учащихся для советской (российской) школы разных поколений, а также сравнить российскую школу с зарубежной}.

Студенты сразу, начиная с первых курсов, получали собственно университетское образование, на старших курсах зачастую делали свою первую научную работу, аспирантура же была сфокусирована не на учебе, а на научной работе.  С явным одобрением об этом упоминал  Л.~С.~Понтрягин \cite{pontriagin}, стр 9--11. И даже сейчас многие студенты мехмата МГУ, при всех современных трудностях, делают свою первую научную работу на 4--5 курсах, чего и близко нет в известных автору западных университетах.

Интерес к научной работе студентов, а также развитие системы внеклассной работы со школьниками породили  вопрос о возможности научной или хотя бы преднаучной работы школьников и её организационных формах. В тридцатые годы 20 века в связи с индустриализацией была  создана система дворцов пионеров, кружков и детского творчества. Примечательно, что решение о
её создании было принято на самом высоком уровне. Конференции школьников постепенно выкристаллизовывались  из
этой системы.

В западных университетах на первых курсах происходит исправление недостатков школьного образования и пополнение необходимых базисных знаний, на старших курсах студенты по существу только начинают осваивать науку как таковую. Срок аспирантуры больше, чем в России, и первая её половина тратится опять же на учёбу, а сама научная работа протекает, как правило, на последних годах аспирантуры и во время ``постдока''. В этой связи для западных студентов приобретают большее значение олимпиады (как школьные, так и студенческие) и оказывается полезен опыт российских конференций школьников.

Вот один только пример. В.~М.~Тихомиров и автор работали (в разное время) в Jacobs University Bremen, где вели курс ``perspectives in mathematics'',
созданный Д.~Шляйхером. Студент оценивался по курсовой работе (term paper).
Оказалось, что деятельность в обоих случаях была чрезвычайно похожа на работу научного руководителя доклада школьника (одна из таких студенческих работ
опубликована в ``Математическом просвещении'', см. \cite{VikolApostol}). То же относится и к курсу, который вел автор в Бар-Иланском университете.

\medskip

Интересно, что вначале кружки мехмата МГУ были основаны на докладах школьников. Современный тип мехматского кружка для школьников является гениальным изобретением Д.~О.~Шклярского и
есть следствие реформы, им проведенной
%\footnote{его и Н.~Н.~Константинова автор считает организационным гением, С.Е.Рукшина -- предгением}
еще до войны. Кружок, основанный на решении задач, стал достаточно эффективным и принес замечательные результаты. Но вместе с тем обучение некоторым сторонам деятельности математика весьма неудобно проводить  в таком формате, так что актуальным стал поиск дополнительных форм работы.

\medskip

{\bf Сложившаяся ситуация.} Сейчас проводится целый ряд научных конференций школьников, как российских так и международных. Активно действует корпорация
Intel ISEF (США): \\ http://www.societyforscience.org/.
Упомянем
европейские конференции школьников: \\ICYS: http://metal.elte.hu/\~icys/, \
EuCYS: http://ec.europa.eu/research/youngscientists/index\_en.cfm. Все большее распространение (в том числе –- по разнарядке ``сверху'') получают ежегодные школьные научные конференции (организуемые, как правило, самими учителями), призванные продемонстрировать успехи в приобщении учащихся к ``исследовательской работе'' в результате внедрения метода проектов. Конференциям посвящены работы  \cite{Plukin,vneshkolnik,Liashko,Skopenkov}.
Следует осмыслить сложившуюся ситуацию, сформулировать
принципы проведения конференций учащихся.\footnote{В 90-е годы автор проводил конференцию ``Поиск'',
проходившую следующим образом.
Автор стремился приглашать сильное и разноплановое жюри. В нем неоднократно принимали участие  Н.~А.~Бобылев, Н.~Б.~Васильев, А.~И.~Галочкин,  А.~А.~Егоров, А.~А.~Заславский,  А.~Г.~Кулаков, В.~Н.~Латышев, Н.~Г.~Мощевитин, В.~В.~Произволов, Г.~Б.~Шабат. Доклады проходили в два дня. Рабочий день состоял из трех частей: а) доклады школьников, б) доклады студентов (по темам, понятным школьникам; цель -- дать представление о реальной научной работе), в) круглый стол --
члены жюри предлагали участникам конференции задачи для самостоятельного исследования.}

%\subsection{Цели автора}

\medskip
{\bf Благодарности.} Автор признателен  А.~К.~Ковальджи, А.~И.~Сгибневу, Б.~Р.~Френкину  за полезные обсуждения. Особая признательность -- Н.~Х.~Розову  за полезные обсуждения и поддержку.

\section{Проектная деятельность и конференции учащихся. Цели и задачи}

%Отношения с учителем
В проектной деятельности отношения с учителем определяются тем, что ученик сам выбирает себе руководителя и участвует в выборе темы. Учитель превращается в консультанта и старшего коллегу, а не просто тренера. В идеале наблюдается сотворчество учащегося и научного руководителя.
При исследовательском подходе к обучению ученик участвует в постановке вопросов, выдвигает гипотезы, которые сам доказывает или опровергает. Полученный ответ может служить основанием для новых вопросов. Результат может быть неизвестен заранее как ученику, так и учителю. В этом особенность  хорошей исследовательской задачи. Работа над ней может идти от нескольких дней до нескольких лет. Очень хорошо, когда используются разнообразные методы и есть простор для частичных продвижений и уточнений, постановки вспомогательных задач.
Невозможно придумать сразу несколько идей одновременно, их надо получить по одной, работая с более простыми ситуациями. К сожалению, спортивный подход, довлеющий над многими олимпиадами \cite{dver}, приводит к тому, что частичные продвижения очень плохо оцениваются.
В результате интеллект развивается, порой весьма значительно, однако при этом человек  разучается мыслить глубоко, со всеми вытекающими последствиями, в том числе академическими.

Методика организации исследовательской работы для начинающих изложена в книге \cite{Sgibnev}, см. также \cite{vneshkolnik,Plukin,Roitberg}. Несколько слов о работе с более продвинутыми людьми (школьниками или студентами). Очень хорошо, если первая работа основана на красивой задаче с понятной мотивировкой, и очень плохо, когда для первой работы руководитель просит подсчитать нечто,
а для чего -- пока непонятно. Первая работа должна оставить ощущение красивого успеха, чего-то праздничного. Трудовые будни -- потом и риск -- потом.
Задача должна быть не слишком простой (неинтересно), но и не слишком сложной -- чтобы человек не сломался. Если все же задача подозрительна на чрезмерную сложность, руководитель должен отслеживать ситуацию и активно помогать. Очень хорошо, когда складываются
деловые отношения -- ученик с учителем вместе хотят получить конкретный результат. Такие отношения обеспечивают подлинное уважение.
В последнее время вышло несколько книг, посвящённых исследовательским задачам для школьников и методике
работы с ними, например, \cite{IvanovRyzik,Origamika,Sgibnev}. Сравнениям конференций и олимпиад посвящена работа \cite{Liashko}.

Задачный подход, при всей своей эффективности, сыграл и отрицательную роль.  Методика проектной деятельности учащихся относительно не разработана, а система -- недостаточно развита. Возможно, при более развитой системе можно получать серьезные научные результаты у старшеклассников значительно чаще,
чем сейчас. В этой связи полезен педагогический опыт Н.~Н.~Лузина. Согласно воспоминаниям Л.~С.~Понтрягина (\cite{pontriagin}, стр. 9--11), педагогический подход Н.~Н.~Лузина состоял в следующем. Первоначальной целью ставилось как можно более быстрое получение научного результата. При этом выбиралась область математики, с одной стороны достаточно значимая, с другой -- не требовавшая обширных знаний. В 20-е --- 30-е годы это была теория множеств, а также теория функций действительного переменного. (В наше время вместо теории множеств можно использовать, например, такую область, как комбинаторика -- в частности,
комбинаторика слов, -- учитывая значительное число бывших олимпиадников.) Такой подход имел очевидные плюсы -- человек, получивший научный результат, приобретал внутреннюю опору. С другой стороны, была опасность замкнуться в узкой области.
\footnote{Человек зачастую не имеет ресурсов на долгую учебу. Очень важно чтобы помимо ``профессиональных'' ученых (т.е. зарабатывающих профессией) были бы и ученые-``любители'' (зарабатывающие на жизнь иначе). Помимо самих научных результатов, это поднимает общественную культуру. Иногда кто-то может прийти в математику и в относительно зрелом возрасте. Им важно быстро получить результаты. Тем самым подход Н.~Н.~Лузина, связанный с относительно быстрым получением первого результата, становится весьма актуальным.
}

Перечислим факторы, создающие потребность в проектной деятельности учащихся.

\begin{enumerate}
    \item Необходимость развивать не только ``пробивные'' возможности, но и исследовательские, творческие качества, креативность.
    \item Есть школьники-``спринтеры'', которые любят быстро решить задачу на олимпиаде и потом отключиться от неё. Есть и школьники-``марафонцы'', которые долго думают над задачей. Для них больше подходят исследовательские задачи (проекты). (Но и для спринтера очень полезно временами бегать марафон.)
    \item Олимпиады плохо приспособлены для ``тугодумов''. Отчасти этим объясняется то, что некоторые специалисты относятся к олимпиадам  скептически. В свое время П.~С.~Александров говорил, что если бы он участвовал в математической олимпиаде, то практически не имел бы шансов получить приз, поскольку может размышлять только в одиночестве, при отсутствии кого-либо рядом, и не способен ``думать на время'', решить задачу к заранее установленному сроку. (См. также п. 7.)
    \item Кроме собственно решения задач, большое
значение имеет еще ``домостроительство'', то есть воспитание теоретического мышления, развитие изобретательности, построение теорий и конструкций.
    \item Очень важна практика в подробном разборе и осмыслении чужих работ, а также
самостоятельном изложении материала. Из попыток ``только просто и доступно'' изложить имеющийся материал часто возникали очень серьезные открытия (например, уравнение Шредингера, схемы Дынкина).
    \item Существенную роль играет возможность приобщения к экспериментальной математике, развивающего эвристическое мышление и воображение.
    \item Хорошо известно, что даже сложная олимпиадная задача становится более доступной при длительном поиске ее решения, особенно при работе в команде.
    \item Научные конференции обеспечивают более плавный переход  от олимпиад к научной деятельности. Отсутствие такого перехода приводит к появлению неудачников среди способных студентов.
    \item На конференциях встречаются ученики, успешные по разным предметам, бывают и комплексные доклады. Конференции являются чрезвычайно полезными в плане межпредметных связей, что важно в дальнейшей деятельности.
    \item Привлечение ученых к работе с молодежью, что не происходит в
достаточной степени на современных олимпиадах \cite{dver}.
В современных олимпиадах произошел разрыв с научным сообществом, происходит процесс вытеснения. Доминируют узкие специалисты по подготовке к соревнованию, не только не являющиеся учеными, но зачастую способные преподавать в весьма специфических условиях \cite{dver}.
\end{enumerate}

\medskip
{\bf Критика.}
Несмотря на вышеперечисленное, можно встретить резко критические отзывы о конференциях учащихся вплоть до утверждений о ``профанации научной работы''. Понимая пристрастность и несправедливость таких утверждений, следует все же выделить в них рациональное зерно.

\medskip
{\bf 1. } Одно из главных критических замечаний  звучит так: `` {\it Очень мало учащихся получают значимые научные результаты.}’’

Однако если обратиться к олимпиадам, то большинство их участников не доходит до высшего уровня, с которым и следовало бы сравнить публикабельные результаты участников конференций. В то же время без основания невозможна вершина: ведь уничтожив олимпиады несложного уровня, мы опустим и высокий уровень. Так же и в научной работе: чтобы прийти к вершинам, необходимы промежуточные результаты, пусть сначала и весьма скромные. Хотелось бы отметить, что перед начинающим школьником не следует сразу ставить особенно
трудную задачу -- она может отпугнуть, быть может, перспективную личность. Для выявления таких личностей среди школьников и учителей важны именно доклады начального уровня -- как и первые этапы олимпиады.

\medskip
{\bf 2.}  К сожалению, современное общество больно, в нем показуха и коррупция проникли во многие сферы жизни, в том числе и в образование, и отнюдь не только в России.
  В погоне за имиджем и во исполнение требований начальства, организаторы олимпиад иногда  ``сливают'' решения задач на олимпиадах (в том числе международных, что неоднократно наблюдал автор), а руководители -- пишут ``доклады’’ за учеников для конференций или предлагают темы, требующие не размышлений, а лишь создания красочных презентаций. Развелось много школьных псевдонаучных конференций, которые только дискредитируют идею, а первые места раздают ``своим’’ (как, впрочем, и жульнических ``олимпиад''). Ситуация усугубляется тем, что далеко не все учащиеся и, увы, учителя имеют адекватное представление о научной этике. Тем не менее представляется, что проблемы показухи и коррупции разрешимы и даже далеко не лучшим образом  организованные конференции, несмотря на ``пену'',  все же приносят определенную пользу, ибо с доклада может начаться интерес к предмету.

 Что есть достойная (пред)научная деятельность учащегося? Здесь возможны варианты.

\begin{enumerate}
    \item Детальный разбор  и самостоятельное переизложение относительно трудного материала (реферативная работа). \footnote{Примером нормального доклада (несложного уровня) служит самостоятельное переизложение теоремы ван-дер-Вардена о раскрасках. Адекватный приз -- брошюра А.~М.~Райгородского, посвященная раскраскам \cite{Raigor}.}

    \item Разного рода компьютерные эксперименты и поисковые лабораторные работы по математике.
\footnote{Темы такого рода лабораторных работ тщательно продумывал А.~Н.~Колмогоров и внедрял их в практику работы школы-интерната при МГУ, которая сегодня носит его имя. Особенно широкие возможности организации лабораторных работ по математике открывает компьютер; см., например, \cite{lab}.}

    \item В домашних (психологически спокойных) условиях (тем более, работая в команде) учащийся способен решить существенно более сложную задачу, чем в условиях олимпиады.
\end{enumerate}

Участники конференций должны проходить не только ``научную’’, но и психологическую подготовку. В частности, следует предостерегать от ожидания немедленных или постоянных успехов. Действующие ученые тоже не каждый день совершают открытия, есть и ``трудовые будни’’ (техническая работа, доведение известных идей и т.п.). Учащиеся зачастую плохо представляют специфику доклада на конференции, а потому чрезвычайно
полезны предварительные консультации (для учеников и учителей) и репетиция доклада.
\medskip

{\bf Предварительные выводы.} Автору представляется, что:
\begin{itemize}
    \item Научные конференции удовлетворяют определенные  потребности и развивают определенные качества учащихся,
не востребованные на олимпиадах.
        \item Эти потребности и качества нуждаются в дальнейшем осмыслении.
    \item Формы организации научной или преднаучной деятельности учащихся не устоялись и требуют
совершенствования.
В целом, однако,
значение конференций пока уступает значению олимпиад. Необходимо мыслить широко, изобретать новые формы работы, экспериментировать. \footnote{В школе всячески надо
пропагандировать привлечение школьников к докладам. Если ученик 8 класса делает доклад
``Золотое сечение'' (чисто реферативный -- пересказывает то, что прочитал и разобрал по
``наводке'' учителя), то это очень хорошо. Он учится самостоятельно читать математическую
литературу (пусть для начала -- популярную), учится ее понимать, учится выступать, слушать
и отвечать на вопросы, объяснять другим непонятное им. Остальные школьники учатся
внимательно слушать, понимать ``со слуха'', задавать вопросы.}

    \item  Однако несовершенство форм не означает отказа от проведения конференций.  Надо прежде всего заботиться о {\it пользе для учащихся}, а уже потом -- {\it о совершенствовании форм проведения}. Дорогу осилит идущий.

    \end{itemize}

\medskip

{\bf Позитивный опыт.} Автору представляется весьма успешным опыт проведения  организацией ``Авангард'' научных конференций школьников
под руководством Д.~В.~Андреева, В.~A.~Тиморина, Е.~Н.~Филатова (конференция
``Интел-Авангард'').
Следует упомянуть
конференцию ``Интел-Юниор'' http://junior.mephi.ru/,
федеральную конференцию ``Юность, наука, культура'', которую организовывал Р.~В.~Плыкин,
конференции лицея ``Вторая школа'' (под руководством П.~В.~Бибикова, А.~К.~Ковальджи, К.~В.~Козеренко), СУНЦ МГУ (``Колмогоровские чтения''), http://www.readings.ru/,
лицея 1502, школы-интерната ``Интеллектуал'', чтения им. Вернадского, конференцию ``Шаг в науку'', организуемую МФТИ,  и др.
Лабораторными  работами и компьютерными экспериментами успешно занимаются школьники на кружке Г.~Б.~Шабата ``Клуб экспериментальной математики''. В 70-е -- 80-е годы (до распада СССР) функционировал замечательный фестиваль науки в г. Батуми, вдохновителем которого была заслуженный педагог Аджарии Медея Жгенти. В работе фестиваля принимали участия члены редколлегии журнала ``Квант'', известные математики. Отчет о работе фестиваля регулярно публиковался в  ``Кванте''.  В 90-е годы конференции школьников проводил замечательный педагог В.~В.~Бронфман. Данный список неполон и автор приносит извинения тем, кого он не включил.

Достаточно совершенной организационной формой служит международная Летняя конференция Турнира городов (ЛКТГ) \cite{lktg}, замечательное изобретение Н.~Н.~Константинова. Школьникам в условиях выездного мероприятия предлагают многошаговые  исследовательские задачи. Это является переходником от олимпиад к научной деятельности. В свое время существовал также заочный конкурс ЛКТГ. К сожалению, он заглох. Одна из причин -- казалось бы, малое число участников. Однако более важной была  подводная часть айсберга -- многие школьники вместе с учителями использовали заочный конкурс для подготовки докладов вне рамок ЛКТГ.
\footnote{Здесь уместно упомянуть удачно найденную И.~М.~Гельфандом форму работы созданной им в 1964 году Всесоюзной заочной математической школы. Эта форма получила название ``коллективный ученик’’ и позволяла учителям отдаленных школ заниматься вместе со своими учениками по программе ВЗМШ и тем самым повышать свою квалификацию.}

Долгое время автор весьма скептически относился к возможности преднаучной деятельности учащихся средних классов школы, тем более не специализированной. Возможно потому, что сам автор сложился в рамках задачного подхода, заданного Д.~О.~Шклярским. Тем не менее работа А.~И.~Сгибнева \cite{Sgibnev} и опыт Красноярской летней школы заставили пересмотреть эту точку зрения.  Была убедительно продемонстрирована возможность и важность проектной деятельности в данной возрастной группе. Педагогические задачи проектной деятельности и ее методика достаточно хорошо изложены в работе \cite{Sgibnev}. Там же довольно много рассказывается и об организационных моментах.

\section{Методические основы проведения конференций}

Некоторые утверждения данной статьи профессиональному математику покажутся очевидными. Однако, с одной стороны, о любом часто встречающемся недостатке надо писать, уделяя ему внимание с целью предупреждения или минимизации. Например, математику очевидно, что докладчик {\it может пользоваться} внешней помощью со стороны учителей, родителей, товарищей, но
обязан это оговорить на докладе. Если же этот принцип не разъяснить, то фактор помощи будет просто скрываться, возникнут неэтичные поступки.  С другой стороны, организатор олимпиад, столкнувшись с таким гораздо менее формализованным мероприятием, как научная конференция, зачастую проявляет профессиональный снобизм (а иногда -- и ревность). Весьма
непроста проблема разношёрстности докладов. Они бывают совершенно разных
типов (см. п. 3.3) и, соответственно, должны по-разному оцениваться (но есть и общие требования). Есть и общие принципы организации работы жюри.

\medskip

\subsection{Принципы работы жюри}

\begin{enumerate}
    \item Ответственность по отношению к учащимся. %  В.~Д.~Арнольд мне как-то говорил, что
    Не надо учащихся ``подставлять’’. Например, сам учащийся может неадекватно оценивать свою работу и ее изложение, так что требование выкладывать работу на сайте www.arxiv.org в качестве предварительного условия окажется провокационным. (Выкладывать работу на таком сайте лучше по решению жюри после конференции.)

    \item Доброжелательность к учащимся. Жюри не должно самоутверждаться за их счет, в частности, проявлять снобизм. Необходима доступность в общении, не следует строить из себя небожителей.
\footnote{Это качество при проведении математических олимпиад всегда проявлял (и требовал того же от своих помощников) А.~Н.~Колмогоров, демонстрируя полную готовность беседовать на любую тему с любым участником.}

    \item Адекватность отбора докладов. Предварительный отсев докладов должен быть {\it максимально мягким}, надо твердо отсекать лишь ``уфологию'' и совсем слабые работы.
Следует доверять известным специалистам и принимать работы по их рекомендации. Относительно слабый доклад может оказаться полезен как участнику, так и его научному руководителю, если они воспользуются советами жюри о направлении дальнейшей работы.
\footnote{Подчеркнем еще раз, что работы начального уровня, подобно первым этапам олимпиады, несут чрезвычайно важную функцию {\it выявления} перспективных личностей -- среди школьников и учителей. Эти люди нуждаются в поощрении за честно выполненные работы, пусть и скромные по меркам высокого профессионала. Кроме того, при организации конференции, особенно при первом ее проведении, чтобы заполучить докладчиков, следует лично обращаться к людям, проявляющим соответствующую активность. Но такие просьбы {\it накладывают определенные обязательства}, в частности {\it необходимость уважать рекомендации}.}

Относительно более слабым докладам (конечно, не сообщая заранее об их ``слабости'')
в интересах слушателей можно выделить меньше времени, что является практикой также ``взрослых'' конференций. Полезно организовать стендовые доклады.
%\footnote{Очень плохой является ситуация
%когда второй диплом получает доклад по уровню не выше чем
%отвергнутый к участию.
%Был случай с участницей из школы-интерната ``Интеллектуал'',
%которая не согласилась с предложенным на консультации способом
%изложения и в результате не была допущена. Это плохо не только тем, что обиделся участник и научный
%руководитель. Главное -- создается (пусть даже ложное)
%впечатление, что консультация вынуждает к приписыванию нового
%научного руководителя.}
% Это соображение, мне представляется, полезно включить!!!

    \item Легитимность жюри. Формирование жюри одним лидером, особенно авторитарным, приводит к разного рода подозрениям. Желательно обеспечить представительство из разных
 организаций, и хорошо когда есть возможность обеспечить представительство из разных городов.

    \item Член жюри (особенно председатель) должен избегать голосования по работам своего подопечного и тем более их не ``пробивать''.

    \item Конференция -- это, прежде всего, {\it праздник науки}.
    {\it Сначала надо искать позитив в содержании работ} (хотя о
недоработках, разумеется, следует сообщать учащимся), {\it а потом уже
оценивать форму изложения.} Учить изложению -– задача научного руководителя, а не жюри. Сегодняшняя школа грамотному и логичному изложению материала учит плохо; распространение тестов усугубляет ситуацию. Школьники (и, кстати, студенты) в этом не виноваты, и к ним тут недопустимо придираться.
\footnote{А.~Н.~Колмогоров и В.~И.~Арнольд помогали своим студентам
написать их первые работы, вплоть до создания текстов. С другой стороны, автору
известен возмутительный случай, когда студент несколько лет
писал пятистраничную статью. Научный руководитель,
вместо оказания помощи, его ругал, тогда как часовая
беседа автора со студентом поправила ситуацию.}

    \item Стиль требований и их уровень у разных людей различаются. Довольно часто
бывает, что те же требования, которые человек предъявляет к себе или своим ученикам,
он распространяет и на участников конференции.  Это иногда весьма неуместно. Следует напоминать жюри, что участники могут быть чрезвычайно чувствительны к оценкам их выступлений и призывать к осторожности.
  %   \footnote{За обращение со своими студентами учёный отвечает результатами своей научно-педагогической деятельности. Если он будет вести себя неправильно, к нему студенты не пойдут. %На конференции такой ответственности нет.}
%\footnote{Кроме того, есть проблема ответственности. Некоторые педагогические идеи могут показаться  разумными, а их ошибочность становится очевидной в поцессе реализации, например, %когда уходят практически все ученики или они вообще покидают математику. Автору известен такой деятель, с весьма агрессивным поведением.}
\footnote{Автору известен случай чрезвычайно активной пропаганды членом жюри своих педагогических идей. Их реализация приводила к тому, что студенты, с которыми он работал, от него уходили или вообще покидали математику. Что касается проведения конференций, их релизация также не привела ни к чему хорошему.}

        \item Оценка учащегося должна определяться {\it позитивными аспектами}: творческими достижениями, красивыми идеями (в том числе методическими), оригинальностью изложения. Не будем забывать, что главная проблема научного сообщества отнюдь не в недостатке пуританской строгости изложения, а в потоке тривиальных результатов, банальных идей и неинтересных статей. Это во многом связано с погоней за числом публикаций и всякого рода сомнительными индексами.

       \item Разумное число докладов. Оптимальное --– примерно 20--25,
а если их меньше 10, то это, скорее, неудача.
\end{enumerate}

%\medskip

\subsection{Организация докладов}

{\bf Общие требования.} Вне зависимости от уровня доклада и его типа, к докладу следует предъявлять следующие требования:

\begin{enumerate}

         \item Доклад должен быть {\it честным}. Конечно, допустимы как переизложение чужих результатов и компиляция, так и смешанные формы с изложением части своих результатов. Допустима помощь родителей, товарищей, учителей и т.д. Но при этом должно быть четко  {\it объявлено, что сделал сам докладчик} и каков вклад третьих лиц.
Докладчик должен также честно изложить
историю вопроса и честно упомянуть персоналии. Если в каком-то месте доклада необходимый анализ исчерпывающе не проведен, это следует прямо отметить, гипотезу нельзя выдавать за доказанный факт, преимущества того или иного подхода надо не голословно объявлять, а убедительно мотивировать.

    \item Докладчик обязан быть {\it квалифицированным}, т.е.
понимать содержание  доклада или хотя бы, как минимум,  той части коллективной
работы, которую он выполнил (в случае, когда доклад {\it комплексный} – включает как математическую, так и программистскую или естественнонаучную часть). Он должен владеть материалом доклада настолько, чтобы отвечать на вопросы по существу изложенного. Вообще говоря, нужно знать даже больше, чем сказано в докладе, чтобы быть готовым к смежным вопросам. Дефекты в понимании содержания доклада служат весьма существенным его изъяном.
\footnote{Бывали ситуации, когда ученик младших классов докладывал ``высшую математику''. Когда обнаруживалось грубое непонимание им
содержания доклада, руководитель школьника атаковал членов жюри, обвиняя их в педагогической некомпетентности, ибо они ``не желают делать поправки на возраст’’.
На это надлежит заметить, что  руководителю следовало  предлагать ученику доступные ему вещи.}

\item Докладчик должен понимать {\it мотивировку},
уметь отвечать на вопросы, почему то, чем он занимается, естественно, почему и как он выбрал данную тему.

\item Докладчик должен иметь личную точку зрения. Содержание доклада должно быть ему не безразлично, действительно глубоко интересно. В противном случае доклад лучше не делать.

     \item Не следует публично решать квадратные уравнения! Доклад -- не урок, где автор демонстрирует свои элементарные технические умения с подробным рассказом у доски. Заинтересованным слушателям можно предъявить более подробные материалы.
\footnote{Это {\it должно быть объявлено} организаторами перед началом конференции. Такое объявление существенно экономит публичное время и улучшает качество
выступлений.}

    \item Докладчик, прежде всего, должен заботиться о качестве изложения,  о понимании
слушателями материала и только затем -- о произведённом на них впечатлении.

\end{enumerate}

\subsection{Типы докладов}
Опыт показывает, что доклады учащихся можно условно разбить на несколько категорий. Каждый тип докладов предполагает свои критерии качества.

\medskip
{\bf 1. Реферативный доклад.} Распространено предубеждение о ``второсортности'' такого рода докладов.
Однако многие крупные научные результаты возникали просто из попыток привести в порядок уже известный материал.
Н.~П.~Долбилин отмечал, что составление хорошего реферата развивает особые качества, тоже важные для математика.
Критерии оценки и требования к такому докладу следующие:

\begin{itemize}
  \item Насколько самостоятельно организован излагаемый материал (а не буквально переписан из книжки); насколько оригинален путь изложения.

  \item Насколько интересна тема.
\end{itemize}

Удачным примером может служить доклад школьника Саши Буфетова по проблеме Варинга, впоследствии опубликованный в журнале ``Фундаментальная и прикладная математика''. % \cite{BufetovKanel}.
Сейчас А.~И.~Буфетов -- д.ф.м.н., профессор НИУ ВШЭ и мехмата МГУ).

\medskip
{\bf 2. Тематический набор задач с решениями.} Пошаговое решение набора мелких задач часто путают с научным исследованием, и потому такого рода доклады могут выдаваться за ``научные''. Школьники зачастую копируют учителей и авторов учебников. В этом случае критерии оценки должны быть
иными.

Следует оценивать оригинальность не только решений, но прежде всего самой подборки задач, объединяющие их идеи.
Если в качестве доклада заявлен задачник, то у жюри возникают вопросы: смотрел ли докладчик с позиций автора задачника на другие книги и учебные пособия, какие выбраны темы и почему.
\footnote{Создание учебного пособия в форме ``многошаговой поднимающейся лестницы’’ является очень удачным приемом, стимулирующим интерес и развивающим креативность обучающегося, но и очень трудным делом. В качестве примера стоит упомянуть классическую книгу \cite{convex}.}
Полезно подчеркнуть, что ученик, преподаватель и автор книги по-разному смотрят на одну и ту же книгу: ученика привлекает доступность изложения материала, преподаватель видит методическую реализацию, а автор оценивает, как написана книга.

\medskip
{\bf 3. Экспериментальная работа.} Обычно она связана с компьютерным моделированием, численным экспериментом и др. В таком докладе оценивается:

\begin{itemize}
  \item Качество постановки эксперимента.

  \item Наличие результатов и их анализ, а также корректность использования статистики. (Заметим: лучше честно признать, что статистические исследования не проводились, чем продемонстрировать грубое непонимание их сути.)

  \item Математическое содержание работы.
  \item Практическая сторона рассмотренной задачи.

  \item Методическая часть (в частности, качество программного интерфейса).
  \end{itemize}

Достоинство работ такого рода -- в относительной их доступности. Возможны комплексные работы, в которых может присутствовать и естественнонаучная часть. %Кроме того, эти
Такие работы легче проводить усилиями целой команды.

\medskip
{\bf 4. Самостоятельное исследование в области чистой математики.} Критерии оценки такой работы -- вкус автора, качество постановки задачи, трудность её решения, новизна полученных результатов. Этот тип работ оценивается наиболее высоко. Важно, чтобы докладчик умел четко объяснять мотивировку и отвечать на смежные вопросы.

Однако полезно иметь в виду следующую возможность манипулирования: учащегося натаскали в некоторой специальной области, указали последовательность утверждений,
которые следуют друг из друга относительно несложным образом, и он все это просто воспроизводит. Поэтому следует выяснить -- понимает ли докладчик {\it мотивировку}.
%\medskip

\medskip
{\bf Оценка докладов.} Конференция -- не олимпиада и не спорт. Ни в коем случае не следует измерять сантиметры и секунды.
\footnote{На конференциях ``Поиск'' и ``Юность, наука, культура'' выдавались два типа диплома -- ``лауреат'' и ``дипломант''. Мы старались не  выдавать более высокие дипломы -- только в исключительных случаях. } Тем не менее важно иметь в виду следующую качественную градуировку:
\begin{itemize}
  \item Доклад, возможно, полезный докладчику и его руководителю. Выдается диплом участника.
  \item Доклад, полезный слушателям. Выдается диплом лауреата.
  \item Доклад, заслуживающий публикации в научном или научно-популярном журнале. Выдается условно первая (иногда вторая -- в зависимости от уровня) премия. Но это только гарнир к основной награде --  публикации.
  \item Получен красивый яркий результат. Выдается первая премия (если таких работ несколько -- то и первых премий несколько, между собой они {\it не сравниваются}). Наградой является помощь при публикации.
\end{itemize}

\section{Заключение. Что и как делать?}

Какие конкретные шаги можно сейчас предпринять?

\begin{itemize}

\item
Надо подумать, как привлекать участников-школьников к действительно научным конференциям, откуда и как они могут появляться.
Это особенно актуально для
отдаленных школ, находящихся вдали от научно-образовательных центров, где и близко нет
никаких квалифицированных руководителей, а уровень учителей часто низок. Как организовать
олимпиаду -- ясно, а как научную конференцию? Сам школьник далеко не всегда может
найти способ с кем-то связаться. Раньше каналом для поиска творческих задач был
журнал  ``Квант''. Надо подумать
над какой-то живой системой заочного консультирования интересующихся школьников
заинтересованными профессионалами. Кто и как может это делать?

    \item Очень полезно собрать и опубликовать список возможных тем. Некоторые темы указаны в книге \cite{Sgibnev}. В данный момент можно предложить организаторам научных конференций учащихся воспользоваться материалами
проектов ЛКТГ (см. http://www.turgor.ru/lktg/ ), а также задачником ``Математического просвещения''. Многие задачи из проектов ЛКТГ до сих пор не
решены, а ряд руководителей проектов согласится на заочные контакты.

    \item Было бы весьма полезно возродить заочный конкурс ЛКТГ.

    \item Представляется целесообразным организовать аналог ЛКТГ для студентов младших курсов и выпускников школ.

    \item  Многие учащиеся (и даже учителя) имеют недостаточное представление о том, как делается доклад на конференции или семинаре, как правильно и доступно рассказывать содержание материала собравшимся слушателям. Докладчика необходимо учить выступать, отвечать на вопросы, дискутировать с оппонентом.  Поэтому полезны консультации докладчиков и тренировки выступления, методическая помощь организаторов конференций учителям и ученикам.

    \item Представляется целесообразным издать своего рода ``темник''~-- сборник возможных тем исследования. Возможна следующая структура книжки. а) Рассказ о конференциях учащихся, критерии оценки докладов. б) Примеры докладов учащихся и студенческих научных статей. в) Темы для рефератов. г) Темы для экспериментальных работ. Задачи, связанные с численным моделированием. д) Исследовательские задачи в области чистой математики. (Здесь полезен опыт ЛКТГ. Задачи должны быть разного типа, в зависимости от уровня подготовки учащихся. Некоторые задачи могут быть близки к школьной программе, другие -- к материалам олимпиад. Следует представить разные темы и разделы математики.) %Частично это реализовано в работе %\cite{Sgibnev}.

\item Попытаться предложить иные формы работы. Данный список неполон, важны дополнительные соображения и идеи. Автор призывает к сотрудничеству.

\end{itemize}

\end{document}